\documentclass{article}
\usepackage{amsmath}
\font\tengoth=eufm10
\font\sevengoth=eufm7
\font\fivegoth=eufm5
\newfam\gothfam
\textfont\gothfam=\tengoth \scriptfont\gothfam=\sevengoth
\scriptscriptfont\gothfam=\fivegoth

\newtheorem{theorem}{Theorem}[section]
\newtheorem{proposition}[theorem]{Proposition}
\newtheorem{definition}[theorem]{Definition}

\newtheorem{lemma}[theorem]{Lemma}

\newtheorem{remark}[theorem]{Remark}

\def\blacksquare{\hbox to .60em{\vrule width .60em height .60em}}

  \font\bb=msbm10 
\catcode`\@=12  

\def\na{\nabla}

\def\?{\'e}
\def\{\`e}
\def\?{\`a}
\def\{\`u}
\def\{\c c}
\def\hb {\hfil \break}
\def\n {\vskip 0.2cm \noindent }
\def\scirc{\,{\raise 0.8pt\hbox{$\scriptstyle\circ$}}\,}
\def\ins{\,{\raise 0.2cm \hbox{ $\scriptstyle \circ$}}\,}

\def  \é{\'e}
\def\è{\`e}
\def\à{\`a}
\def\ù{\`u}
\def\ç{\c c$\!\!\!$}

\date{ }
\title{ Global structure of webs in codimension one   
   }
\author{Vincent Cavalier \and Daniel Lehmann}
 
\begin{document}

\centerline{\bf   GLOBAL STRUCTURE OF WEBS IN CODIMENSION ONE  }
 
 \bigskip
 
  \rightline{(\état  du $17/03/2008)$}
  
  \bigskip
  
 \centerline{ Vincent Cavalier and Daniel Lehmann}
 
 \vskip 2cm
 
  \begin{abstract}
  
  In this paper, we     describe the global structure of webs in codimension one, and study in particular   their singularities (the  \emph{caustic}). We define and study different concepts which have no interest locally near regular points, such as the \emph{type}, the \emph{reducibility}, the \emph{quasi-smoothness}, the \emph{CI-property} (complete intersection), the \emph{dicriticity}....   As an   example  of  application, we explain 
 how the algebraicity of a global holomorphic web in codimension one on the complex $n$-dimensional projective space {\bb P}$_n$ (algebraicity that we prove by the way to be equivalent to the linearity) depends only on the behaviour  of this web  near its   caustic, at least for quasi-smooth  webs with CI irreducible components.
  
 \end{abstract}
 
  \vskip 3cm

 1- Introduction 
 
  \bigskip
 
 2- Background on the manifold of contact elements 
 
  \bigskip
 3-  Webs of codimension one on a holomorphic manifold $M$ 
 
 \bigskip
   
  4- Global polynomial partial differential equations

   \bigskip
   
   5- Complete and locally complete  intersection webs  (CI and LCI)
   
     \bigskip
     
    6- Webs on the complex projective space {\bb P}$_n$ 
     
       \bigskip
       
       7- Dicriticity    and algebraicity of global  webs on {\bb P}$_n$ 
       
         \bigskip
         
         References
         
         \vfill
         
   {\n A.M.S. Classification : 57 R 20, 57 R 25, 19 E 20.\hb
 }

\maketitle
\section{Introduction}

In this paper, we want   for instance  to explain 
 how the algebraicity of a global holomorphic web in codimension one on the complex $n$-dimensional projective space {\bb P}$_n$ reflects  on the behaviour  of this web  near its singularities, and -in some cases- can be readen on this behaviour, generalizing a theorem already proved in [CL1] for $n=2$.  
 By the way, we shall need to describe the global structure of webs in codimension one, which may be interesting in itself. We shall define in particular   the \emph{quasi-smoothness}   and the \emph{dicriticity}, both concepts   which can be  readen on the singular part   of the web (its \emph{caustic}). In fact, we   shall see easily that any linear web on {\bb P}$_n$  is dicritical (in fact, any globally defined linear web is algebraic).    Conversely, at least for webs whose any irreducible component is  a \emph{complete intersection} (CI),  we shall prove that  the quasi-smoothness (a   generically satisfied condition) and the  dicriticity (generically not satisfied) imply    the algebraicity. The interest of such a result is to provide a situation where properties of the web on the caustic imply properties everywhere else. 
 
  The general method will be to define a web of codimension one on a $n$-dimensional manifold $M$ by the data of a $n$-dimensional subvariety $W$ of the manifold $\widetilde{M}$ of the contact elements of $M$, the tautological contact form on $\widetilde{M}$ being integrable on the regular part $W'$ of $W$, and    inducing consequently a foliation $\widetilde{\cal F}$ of codimension one.   This foliation  is nothing else but a kind of desingularization (``decrossing") of the web:   above the regular part $M_0$ of the web (the points of $M$ where the leaves  of the local foliations defined by the web are not tangent),the  leaves of $\widetilde{\cal F}$ project locally onto the leaves of the local foliations in $M$, the part $W_0$ of $W$ above $M_0$ being simply a covering space of $M_0$. But we shall be also interested by looking what happens above the part $\Gamma_W$ of $W $ which projects onto the caustic $M\setminus M_0$.  The quasi-smoothness means that any irreducible component $C$ of $W$ is smooth ($C'=C$), so that we may use on $C$ the general tools for foliations on smooth manifold, and in particular the vanishing theorem of Bott : dicriticity means in fact that the foliation $\widetilde{\cal F}|_C $ induced by $\widetilde{\cal F}$   on $C'$ has no singularity ; in this  case,  the square $(c_1)^2 \bigl(N(\widetilde{\cal F}\vert_C)\bigr)$ of the Chern class of  the    bundle $N({\widetilde{\cal F}\vert_C})$  normal to the foliation ${\widetilde{\cal F}\vert_C}$, must vanish. In particular, when $M=$ {\bb P}$_n$, this will imply that all partial degrees $ \delta_\alpha $  of the CI-web are necessarily zero : this characterizes algebraic CI-webs.

  Let us   be now more precise. Let $d$ be an integer $\geq 1$. Locally,  on  an open set  $U$  of {\bb C}$^n$ with coordinates $x=(x_\lambda)_\lambda$, a codimensional one holomorphic $d$-web on   $U$   is a family of $d$     one codimensional holomorphic foliations ${\cal F}_i$, $(1\leq i\leq d)$, on $U$.  The regular part  of the web is the open subset $U_0$ of $U$ where these  foliations are  all non-singular and where their leaves are not tangent. We shall not require in this definition of $U_0$   that any $k$ of these $d$ foliations  are in general position  for  $k\leq n $. Each  foliation ${\cal F}_i$ may be given  by  a holomorphic 1-form $\omega_i=\sum_{\lambda=1}^n a_i^\lambda(x)\ dx_\lambda$ which is integrable ($\omega_i\wedge d\omega_i\equiv 0$), defined up to multiplication by a unit (holomorphic nowhere-vanishing function $u$), and whose kernel is the Lie algebra of  vector fields tangent to ${\cal F}_i$.  If we assume $a_i^n(x)\neq 0$, the leaves of ${\cal F}_i$ will be the hypersurfaces which are the graphs $x_n=\varphi(x_1,\cdots, x_{n-1})$  of  the common solutions $\varphi$ to  the $n-1$ partial differential equations  $a_i^\alpha (x)+a_i^n(x){{\partial \varphi }\over{{\partial x_ \alpha}}}=0$, $(1\leq \alpha\leq n-1)$. 
   
   We can still say that a local $d$-web is defined  
  by the union $W= \bigcup_{i=1}^d U_i$   of the $n$-dimensional   $U_i$'s, each   $U_i$ being  defined by the $n-1$ equations  $a_i^\alpha (x)+a_i^n(x)p_\alpha=0$ ($1\leq\alpha\leq n-1$)  in $U\times  \hbox{\bb C}^{n-1}$ with coordinates $\bigl(x=(x_\lambda),p=(p_\alpha)\bigr)$,   ($1\leq\lambda \leq n\ ,\ 1\leq\alpha\leq n-1$). Observe that the local contact form $dx_n-\sum_\alpha p_\alpha\ dx_\alpha$ has a restriction to $W$ which is integrable, and that the foliation  $\widetilde{\cal F}$ defined by it on $W$ has a restriction to $U_i$ which projects onto $ {\cal F}_i$ by the map $(x,p)\mapsto x$. Of course, this transcription in $U\times  \hbox{\bb C}^{n-1}$    may seem   pedantic and needless while  we remain locally on $U$. But when
  we wish to define a $d$-web globally on a holomorphic $n$-dimensional manifold $M$ which is no more necessarily equal to an open set $U$ of {\bb C}$^n$, it may happen, as already   observed in the case $n=2$  (see [CL1]),  that the local foliations   ${\cal F}_i$ are not    globally distinguishable. Thus, it will be useful to handle them all together,    $(x,p)$   becoming  local coordinates  on the manifold $\widetilde{M}$ of contact elements of $M$, and $W$ a $n$-dimensional subvariety of $\widetilde{M}$, on the smooth part $W'$ of which the tautological contact form of $\widetilde{M}$ becomes integrable,   $W_0\to M_0$ becoming  a $d$-fold covering space  which is  no more necessarily trivial.  
   
 In the next section 2, we recall some basic facts, more or less known, on the $2n-1$-dimensional manifold $\widetilde M$ of contact elements in $M$.  
 
In section 3, we give the general  definition of a  global $d$-web on $M$ as    a $n$ dimensional subvariety $W$ of $\widetilde M$, globalizing the previous considerations, and precising the  \emph{critical set} $\Gamma_W$ and its projection onto $M$ which is the  caustic.   Various concepts such as the  \emph{type} of a $d$-web, the reducibility (decomposition of  the   web as the juxtaposition of a finite number of global ``irreducible" webs), the ``smoothness"  and the ``quasi-smoothness", the ``dicriticity", etc... are defined. 
 
 After defining   a global PDE in section 4 as a particular  hypersurface $S$ in $\widetilde M$, we study in section 5 the  particular case of the so-called \emph{CI-webs}  (when $W$ is the complete intersection $W=\cap_{\alpha=1}^{n-1}S_\alpha$ of $n-1$ distinct global PDE's $S_\alpha$), and \emph{LCI-webs} (locally complete intersections). Many calculations are much easier in these cases. In particular, after   defining  the \emph{critical scheme}, whose underlying analytical set is $\Gamma_W$, we give a practical criterium in this case for a web to be dicritical. 
 
In section 6, we study some properties  specific to the case $M=$ {\bb P}$_n$. Let  {\bb P}$'_n$ be the dual projective space of hyperplanes in {\bb P}$_n$.  Both $\widetilde{\hbox {\bb P} _n }$ and  
$\widetilde{\hbox {\bb P}'_n }$ are equal to the sub-manifold  of $\hbox {\bb P} _n \times \hbox {\bb P} '_n $ whose points are the pairs $(m,H)$ of a point $m\in \hbox {\bb P} _n $ and a hyperplane $H$ such that $m\in H$. An  interesting duality exists between    webs on {\bb P}$_n$ and    on  {\bb P}$'_n$. In particular,  
some sub-variety $W$'s of $\widetilde{\hbox {\bb P} _n }$ may define  simultaneously  a non-dicritical web on {\bb P}$_n$ and a non-dicritical web on {\bb P}$'_n$ (the so called \emph{bi-webs}).  We shall study also the degree of a web  and the multi-degree of a CI-web: in particular, the CI-webs which are algebraic are those of muti-degree 0.
Their  leaves  are also the hyperplanes tangent to an algebraic developable hypersurface (hypersurface equal to the union of a one-parameter family of $(n-2)$-dimensional projective subspaces of {\bb P}$_n$, along which the tangent hyperplane remains the same).
The \emph{linear} webs are those whose all leaves are hyperplanes of {\bb P}$_n$ : they all are dicritical. 
The cohomology of   $\widetilde{\hbox {\bb P}_n }$   is then computed. 

We have now all the²² needed tools for giving   in section 7 the properties of the dicritical webs on {\bb P}$_n$. 
  
 In a second part of our study which will be published somewhere else (see preprint [CL2]), we   focus on  abelian relations. We   study   there   the webs satisfying to a   condition which is generically satisfied, the   so called \emph{regularity}:
 
 -  the rank of such a  regular web at a non-singular  point (i.e. the dimension   of the vector space of germs of abelian relations    at that  point) is upper-bounded by some  number $\pi'(n,d)$  which, for $n\geq 3$, is strictly smaller than the number  $\pi(n,d)$ of Castelnuovo (the maximal geometrical genus of an algebraic  irreducible and non-degenerate curve of degree $d$ in 
{\bb P}$_n$,  which    Chern proved to be the upper-bound  of the rank in the general case ([C]));

- we define on   $M_0$, when $d$ is equal   for some $r$ to the the dimension $c(n,r)$   of the  vector space of all homogeneous polynomials of degree $r$ in $n$ variables with scalar coefficients,  a generalization of the Blaschke curvature; this is the curvature of some connection. The non-vanishing of this curvature is an obstruction for the rank of the web to reach the maximal value $\pi'(n,d)$. 

\n [In the case $n=2$, any web is regular, any $d$ may be written $c(2,d-1)$, the numbers $\pi (2,d)$ and $\pi'(2,d)$ are equal;  we then recover, globalizing them on $M_0$,  the Blaschke curvature defined  locally in [H] for $n=2,\ d\geq 3$, generalizing the case $n=2, d=3$ of Blaschke ([B])].  
  
\section{Background on the manifold of contact elements}
 
Let  $M$ be a  holomorphic (non-singular) manifold, not necessarily  compact. Denote by $n$ its  complex dimension, $TM$  its complex    tangent bundle, and
 $\widetilde M $ the $2n-1$-dimensional manifold equal to the  total space of the grassmanian bundle   {\bb  G}$_{n-1} M\buildrel \pi\over\rightarrow M$  (a point of $\widetilde M $  over a point $m\in M$  is a  $n-1$-dimensional sub-vector space   of $T_mM$, and is  called  {\it  a contact element} of  $M$ at $m=\pi (\widetilde m)$. 
 
 In the sequel, indices such as $\lambda,\mu,\cdots$  will denote integers runing from $1$ to $n$, while   $\alpha,\beta,\cdots$ will denote integers runing from $1$ to $n-1$.
 
  For all family $v=(v_ \alpha)$ of $n-1$ vectors linearly independant in   $T_mM$, $[v]$ or $[v_1,\cdots, v_{n-1}]$ or  $[(v_ \alpha)_{_ \alpha}]$   will denote the contact element  generated by $v$  in    $\widetilde M_m  $.
 
  Since any $\widetilde m$ over $m$ is the kernel of a non-vanishing 1-form on $T_mM$ well defined up to multiplication by a scalar unit, we can also identify $\widetilde M $ with
the total space of  {\bb  P}$T^*M\to M$ (the projectivisation   of the complex cotangent bundle $T^*M$).

\subsection{Canonical bundles and tautological contact form on $\widetilde M$   :}

  Let
 ${\cal T} \subset  \pi^{-1}(TM)$ be the tautological vector-bundle over $\widetilde M $ (whose fiber over  $[v_1,\cdots, v_{n-1}]$ is the sub-vector space of $T_mM$ generated by ($  v_1,\cdots, v_{n-1} $), identified  with the line in $T_m^*M$ of all 1-forms vanishing on the $v_\alpha$'s ($1\leq \alpha\leq n-1$). 
\begin{lemma} The dual ${\cal L}^*$ of the quotient bundle ${\cal L}=\pi^{-1}(TM)/{\cal  T}$ is the  tautological complex line-bundle of {\bb P}$( T^*M)$.
 \end{lemma} 
  The proof is obvious.

Let $(*)$ and $(**)$ the two obvious  exact sequences of vector bundles
$$(*){\hskip 2cm}0\to  {\cal T}\to \pi^{-1}(TM)\to  {\cal L} \to 0  .$$
 and $$(**){\hskip 2cm}0\to {\cal V}\to T\widetilde M\to \pi^{-1}(TM) \to 0 ,$$
where $\cal V$ denotes the sub-bundle of tangent  vectors to $\widetilde M $ which are "vertical" (i.e. tangent to the fiber of $\pi$).

By composition of  the projection $T\widetilde M\to \pi^{-1}(TM)$ of $(**)$ with the projection $\pi^{-1}(TM)\to  {\cal L}$ of $(*)$, we get a canonical   holomorphic 1-form  on $\widetilde M$, with  coefficients in  ${\cal L}$ $$\omega:T\widetilde M\to{\cal L},  $$   which is called  {\it the tautological contact form}. This terminology will be justified below by lemma 2-3.
\begin{remark}
 We   imbed  ${\cal L}^*$ into $\pi^{-1}(T^*M)$ by duality of $(*)$, and $\pi^{-1}(T^*M)$ into $T^*\widetilde M$ by duality of $(**)$. We may also imbed  ${\cal T}$ into $T\widetilde M$, when we identify it to the kernel of $\omega$. 
 \end{remark}

\subsection{Local coordinates  on    $\widetilde M$}

 Let  $x=(x_1,x_2,\cdots x_n)$ be  local  holomorphic   coordinates on an   open set  $U$ of  $M$, and $m$ a point of $U$. Let  $u= (u_1,\cdots,u_n)$   denote the coordinates in $T_m^*M$ with respect to the basis $(dx_1)_m, \cdots,(dx_n)_m$;  let  $ [u]=[u_1,\cdots,u_n]$ denote the point in {\bb P}$T_m^*M$ of     homogeneous coordinates $u$
 with repect to the same basis  in $T_m^*M$,   and  let  $p=(p_1,p_2,\cdots p_{n-1})$ be the system of affine coordinates on the affine subspace $u_n\neq 0$ of {\bb P}$T_m^*M$ defined by $p_\alpha=-{{u_\alpha}\over{u_n}}$ \   ($1\leq \alpha\leq n-1$).  
 
 Therefore, $(x,p)=(x_1,x_2,\cdots x_n, p_1,p_2,\cdots p_{n-1})$ is a system of local holomorphic  coordinates on the open set $\widetilde {U}$ of contact elements above $U$ which are not parallel to $ \bigl({\partial\over{\partial x_n}} \bigr)   $ : the point of coordinates $(x,p)$ is the $(n-1)$-dimensional subspace of $T_m^*M$ which is the kernel   of the 1-form    $\eta=(dx_n)_m-\sum_{\alpha=1}^{n-1}p_\alpha(dx_\alpha)_m$ ; it is  generated by  the 
the vectors   $(X_\alpha)_m= \Bigl({\partial\over{\partial x_\alpha }}\Bigr)_m  +  p_\alpha  \Bigl({ \partial\over{\partial x_n}}\Bigr)_m$,\break $(1\leq  \alpha\leq n-1)$.

Obviously, the vector fields $X_\alpha= \Bigl({\partial\over{\partial x_ \alpha}}\Bigr)  +  p_ \alpha \Bigl({ \partial\over{\partial x_n}}\Bigr)$, $(1\leq  \alpha\leq n-1)$, 
define a holomorphic local trivialization $\sigma_{\cal T}$ of ${\cal T}$ above $\widetilde {U}$, while the image $\sigma_{\cal L}=\bigl[  {\partial\over{\partial x_n}} \bigr]$ of $  {\partial\over{\partial x_n}} $ by the projection $\pi^{-1}(TM)\to  {\cal L}$ of $(*)$ defines a holomorphic local trivialization   of ${\cal L}$.  
\begin{remark}

By permutation of the order of the coordinates $x_i$, such that the new last one be no more the former last one $($writing for example $  x'_n=x_1,\ x'_1=x_n$ and  $x'_ \lambda=x_ \lambda $ for $ \lambda\neq 1,n )$, we   get a new system of local coordinates on some open set $\widetilde {U'}$ containing contact elements parallel to $ \bigl({\partial\over{\partial x_n}} \bigr)_{m }   $.  
 
 \end{remark}   
      
\begin{lemma}{\it

\item $(i)$ The local   form $\eta=dx_n-\sum_{ \alpha=1}^{n-1}\ p_ \alpha\ dx_\alpha   $    is a contact form  on $\widetilde {U}$, and defines a local holomorphic trivialization of ${\cal L}^*$. Moreover,  the trivializations $\eta$ and $\sigma_{\cal L}=[  {\partial\over{\partial x_n}} ]$
are dual to each other.

\item $(ii)$ The  1-form $\eta\otimes \sigma_{\cal L}$ is equal to the restriction of $\omega$ to $\widetilde U$.} 
\end{lemma}

\n  {Proof :}
 Both forms   $\eta\otimes \sigma_{\cal L}$ and $\omega  $ vanish  when applied to the vertical vectors and to the vectors $   {\partial\over{\partial x_ \alpha}}   + p_\alpha   \ {\partial\over{\partial x_n}}$    (with $ \alpha\leq n-1   $), and take the value $\sigma_{\cal L}$ when applied to ${\partial\over{\partial x_ n}}$,   hence are equal. Moreover  $\eta$ belongs to the dual 
 ${\cal L}^*$ (identified to the set of 1-forms on $\widetilde {M }$ vanishing on $\cal T$).  Since $\eta({\partial\over{\partial x_ n}})\equiv 1$, $\sigma_{\cal L}$ and $\eta$ are dual to each other. 
We check easily that $\eta\wedge (d\eta)^{n-1}$ is a volume form.

\rightline {QED}

\subsection{Change of local coordinates and transition functions}

Let  $(x', p')$ (or $(x',u')$) be  the local coordinates   associated to  $x'=(x'_1,\cdots,x'_n )$
 above some open set  $U'$ of $M$ such that $U\cap U'\neq \emptyset$. Let  $J$ (or $J(x,x') $ in case of ambiguity)    be  the jacobian matrix $J={{D(x'_1,\cdots, x'_n) }\over {D(x _1,\cdots, x _n)}} $ with coefficients $J^\lambda_\mu={{\partial  x'_\lambda  }\over {\partial  x_ \mu  }}$. Denote by  $ {\Delta_n^M}$ (or 
  $ {\Delta_n^M}(x,x')$) its determinant, and let \break $K$ (or $K(x,x')$ be the matrix $K= {\Delta_n^M}\ .\ J^{-1}$, with  coefficients 
    $K^\lambda_\mu={\Delta_n^M} . {{\partial  x_\lambda  }\over {\partial  x'_ \mu  }}$.  In the sequel $u=(u_1,u_2, \cdots,u_n)$ will also be  understood as a line-matrix $(u_1\  u_2\ \cdots\ u_n)$.

 \begin{proposition}
 \hb 
 $(i)$ The coordinates $u$ and $u' $ are related by the formulae $u'=( 1/{\Delta_n^M})\ u\scirc K$ and $u =    u' \scirc J )$, while for  homogeneous coordinates $[u']$ $($resp. $[u ])$, defined up to multiplication by a non-zero scalar,  we may write  $[u']= [ u \scirc K]$ 
 $($resp. $[u] =    [u' \scirc J] )$. The following formulae   hold  : 
   $$p'_\alpha=- \ {{  \sum_{ \beta }  K^\alpha_\beta \  p_\beta-K _n^ \alpha  }\over{{ \sum_{ \beta }  K_ \beta ^n    \ p_ \beta-K_n^n   }}}  \hbox{ and } p _\alpha=- \ {{  \sum_{ \beta }  J^  \alpha_ \beta\  p'_\beta-J_n^\alpha  }\over{{ \sum_{ \beta }  J_\beta ^ n    \ p'_ \beta-J_n^n   }}}. $$

 \hb $(ii)$  The expressions $\Delta_{n-1}^ {\cal T}(x,x') =-\bigl( \sum_{ \beta }  K_ \beta ^ n     \ p_ \beta-K_n^n  \bigr) $   constitute     a system of  transition functions for the bundle $\bigwedge ^{n-1} {\cal T}$ :
 $$X_ {1}\wedge\cdots \wedge  X_ {n-1}=\Delta_{n-1}^ {\cal T}(x,x')\ \bigl(X'_ {1}\wedge\cdots \wedge  X'_ {n-1}\bigr).$$

  \hb $(iii)$ Let $\eta'=dx'_n-\sum_{ \alpha} \ p' _\alpha\ dx'_ \alpha$.  The following formula holds : 
    $$\eta'=\bigl({{\Delta_n^M}(x,x') /{\Delta_{n-1}^ {\cal T}(x,x')}}\bigr) \ \ \eta ,$$   and   the functions ${{\Delta_n^M}  /{\Delta_{n-1}^ {\cal T} }}$ constitute a system of transition functions for ${\cal L}$ : 
   
    \end{proposition} 
    
   \n Proof :

   In fact the formula  $u =    u' \scirc J$ is equivalent to the definition of the jacobian matrix. Interchanging $x$ and $x'$, we get $u' =    (u  \scirc J^{-1})$ ; but, at the level of the projective space $J^{-1}$ and $K$ induce the same automorphism,  and   $[u']$ is also equal to $  [ u  \scirc K]$.
   
  We then deduce the   following   formula, by writing  $u_n=u'_n=-1$, $u_\alpha=p_\alpha$ and $u'_\alpha=p'_\alpha$.
Interchanging $J$ and $J^{-1}$, we get also the $p_\alpha$'s in function of the $p'_\alpha$'s.  

Let $\Delta_{n-1}^ {\cal T}$ be a transition function for the bundle $\bigwedge^{n-1} {\cal T}$. From the exact sequence 
$(*)$, we deduce the isomorphism 
     ${\cal L}\cong \bigwedge^{n}\pi^{-1}(TM)\otimes\bigwedge^{n-1} {\cal T}^*    $. Therefore, ${{\Delta_n^M} /{\Delta_{n-1}^ {\cal T}}}$ is a system of transition functions for ${\cal L}$. 
     
     Moreover, 
 since both 1-forms $\eta$ and $\eta'$ are colinear, we deduce the equality $$\eta =   {{- 1} \over  {\Delta_n^M}}\   \Bigl({ \sum_{ \beta }  K_ \beta ^ n     \ p_ \beta-K_n^n  } \Bigr)  \ \eta'  $$ from the equality $\Delta_n^M.dx_\lambda=\sum _\mu  K^ \lambda_ \mu\ dx'_\mu$. This proves that $\Delta_{n-1}^ {\cal T}=
 -\bigl( \sum_{ \beta }  K_ \beta ^n     \ p_ \beta-K_n^n  \bigr) $ and achieves the proof of the proposition.

      \section {Webs of codimension one on a holomorphic manifold $M$} 
      
Let $W$ be a sub-variety of $\widetilde M$ having  pure  dimension $n$. Let $\pi_{_W}:W\to M$ be the restriction  to $W$  of the projection $\pi : \widetilde M\to M$ . Denote by

$W'$  the regular part of $W$,

$\Gamma_{_W}$ the set of points $\widetilde m\in W$ where
\hb  ${\hskip 1cm}$- either $W$ is singular,
 \hb ${\hskip 1cm}$- or the differential $d\pi_{_W}:T_{\widetilde m}W'\to T_{\pi(\widetilde m)}M$  is not an isomorphism.

 $W_0$   the complementary part (included into $W'$) of $\Gamma_{_W}$ in  $W $,

 $M_0=\pi(W_0)$.
 
 Let $d$ be an integer $\geq 1$.

\begin{definition}

 We shall say that $W$ is a \emph{$d$-web} if 
 
 $(o)$ The map $\pi_W:W\to M$ is surjective \footnote{This condition is automatic when $M$ is compact, because of $(ii)$ and $(iii)$.}.

$(i)$ The restriction $\omega_{_W}$ to $W'$ of the tautological contact form $\omega:T\widetilde M\to {\cal L}$   satisfies to the integrability condition:  this means that, given a local trivialization $\sigma_{\cal L}$ of $\cal L$, the restriction $\eta_{_{W }}$ to $W'\cap \widetilde{U}   $     of the  local  contact form  $\eta$ such that $\omega=\eta\otimes\sigma_{\cal L}$  satisfies to the condition $\eta _{_W}\wedge d \eta_{_W}\equiv 0  $ $($this condition not depending   on the local trivialization   $\sigma_{\cal L}$$)$.

$(ii)$ The restriction $\pi_{_W}:W_0\to M_0$ of $\pi_{_W}$ to $W_0$ is a $d$-fold covering 
$[$The integer $d$ will be called the  {\bf weight} of the web$]$. 

$(iii)$ The analytical set $\Gamma_W$ has complex dimension at most $n-1$, or is empty.

 $(iv)$ For any $m\in M$, the set $(\pi_{_W})^{-1}(m) =W\cap (\pi) ^{-1}(m)$ is an algebraic subset of degree $d$ and dimension 0 in 
 {\bb P}$T^* _m(M)$. 
 
  The analytical set $\Gamma_W$ is called \emph {the critical set} of the web, its projection $\pi(\Gamma_W)$ \emph{the caustic} or \emph{the singular part}, and is complementary part $M_0$ \emph{the regular part} of the web.  
 \end{definition}
\begin{remark} \hb
  $(i)$ Any    $d$-web $W$ on $M$ induces a $d$-web  $W|_U=W\cap \pi^{-1}(U)$ on any
open set $U$ of $M$.
\hb $(ii)$  In this definition, we do not require for $d\geq n$ that any $n$-uple of
points in a same fiber of \break $W_0\to M_0$ be in general position : it may exist
some point $m\in M_0$ and some $n$-uple  $(\widetilde m_1,\cdots, \widetilde m_{n})$
of points  in $(\pi_{_W})^{-1}(m) $ which are in a same  $(n-2)$-dimensional
sub-vector space  of
 the $(n-1)$-dimensional vector space $ \bigl(\hbox {\bb  G}_{n-1} M\bigr)_m $.

 \end{remark}
 \begin{definition}
 The irreducible components $C$ of the variety $W$ in $\widetilde M$ are called the \emph{components} of the web. They all are webs.   The web is said to be \emph{irreducible} if it has only one component. It is said to be \emph{smooth} if $W$ is smooth $(W'=W$, which implies that it is irreducible$)$, and \emph{quasi-smooth} if any irreducible   component  is smooth.
 \end{definition}

    \subsection{The foliation  $\widetilde{\cal F}$ }
  The integrability condition  (i) of the definition of a web   implies that the distribution   of vector fields on $W'$ belonging to the kernel of $\omega_{_W}:TW'\to {\cal L}$ is   involutive: the integral submanifolds of this distribution defines therefore a foliation $\widetilde{\cal F}$ on $W'$. This foliation may have singularities, but not on $W_0$ : in fact, the restrictions $\widetilde x_i$ of $x_i$ to $W_0$ ($1\leq i\leq n$) define local coordinates on every sheet of the covering  $W_0\to M_0$.  With respect to these coordinates,  $\eta_{_{W_0}} $ is written $d\widetilde x_n-\sum_{\alpha =1}^{n-1}p_ \alpha(\widetilde x)\ d\widetilde x_ \alpha$ and does not vanish.

  By projection $\pi_W$ of the restriction   of this foliation to $W_0$, we get locally, near any point $m\in M_0$, $d$ distinct one-codimensional regular  foliations $ {\cal F}_i$, and $\widetilde{\cal F}$ may be understood as a ``decrossing'' of these $d$ foliations. 
  \begin{definition} We call \emph {leaf of the web} any hypersurface in $M $ whose intersection with $M_0$ is locally a leaf of one of the local foliations $ {\cal F}_i$.
   \end{definition} 
  
  Conversely, given  $d$ distinct regular foliations $ {\cal F}_i$ ($1\leq i\leq d$) of codimension one   on some open set  $U$ of   $M$, assume    moreover that these foliations are mutually transverse to each over at any point of $U$, and that there exists    holomorphic coordinates $(x_1, \cdots, x_n)$ on $U$ such that ${\partial\over {\partial x_n}}$ be tangent to none of the foliations $\widetilde{\cal F}_i$:  we may define   $ {\cal F}_i$ by an integrable non-vanishing 1-form $ \eta_i=d  x_n-\sum_{\alpha =1}^{n-1}p^i_ \alpha(  x)\ d  x_ \alpha$.   The sub-manifold  $U_i$ of $\widetilde M$ defined by the $n-1$ equations $p_\alpha=p^i_ \alpha(  x)$ ($1\leq  \alpha\leq n-1$) does not depend on the choice of the local coordinates. The union $W_U=\coprod_{i=1}^dU_i$ is then a $d$-fold trivial covering space of $U$, and defines  a $d$-web (everywhere regular, with no critical set). Denoting by $\pi_i:U_i\to U$ the restriction of $\pi$ to the sheet $U_i$,  the form $(\pi_i)^*(\eta _i)$ is equal to the restriction of $d\widetilde x_n-\sum_{\alpha =1}^{n-1}p_ \alpha(\widetilde x)\ d\widetilde x_ \alpha$ to $U_i$. 
  
  \subsection{Dicriticity : } 
 
   Since $TW_0$ is naturally isomorphic to   $(\pi_W)^{-1}(TM_0)$,  the natural injection ${\cal T}   \to (\pi)^{-1}(TM)$  defines   an obvious injective map ${\cal T}  |_{W_0}\to TW_0$, and the image of this map is the tangent bundle to the foliation $\widetilde{\cal F}|_{W_0}$ since it is annihilated by $\omega_W$ : the foliation 
   $\widetilde{\cal F}$ is non-singular 
  on $W_0$. 
   
    Let now $C$ be a irreducible  component of $W$ and $C'$ its regular part. In general,  $\widetilde{\cal F}$ will have  singularities on $C'$.  
    \begin{definition} We shall say that the  web is \emph{dicritical on $C$}
   if the restriction of the foliation $\widetilde{\cal F}$ to $C'$ is non-singular. We shall say that the web is \emph{dicritical} if it is dicritical on every  irreducible component. 
      \end{definition}
      
      The dicriticity on $C$ is equivalent to the possibility of extending ${\cal T}  |_{W_0}\to TW_0$ 
      as an injective  morphism 
   ${\cal T}  |_{C'}\to TC'$.  We shall give below a practical criterium for that, in the case of LCI-webs. 
      
      We shall see also that any totally linearizable   web on a manifold $M$ (for example any linear  web on {\bb P}$_n$) is dicritical.

  \subsection{Irreductibility of webs and indistinguishability of the foliations : }

The results of this subsection will not be used below. 

Let's consider 
  
  - the number $N_1$ ($\geq 1$) of irreducible components of  $W$,
  
  - the number $N_2$ ($\geq 1$) of connected components of $W_0$, 
  
  - the number $N_3$ ($\geq 0$) of global distinct foliations on $M_0$ whose leaves coincide locally near a point with a leaf of one of the
  $d$ local foliations ${\cal F}_i$ near this point. 
  
 It is clear that $N_3\leq N_2$. [We may have  $N_3<N_2$ ; consider for instance,  the 2-web  on $M=$\ {\bb P}$_2$ whose leaves are the tangent lines to a given  proper conic $X$. Then, the caustic is $X$, and the   2-fold covering $W_0\to M_0$ is not trivial ;   however, $N_3=0$, while $N_1=N_2=1$]. 
 
 In particular, when
  $M_0$ is connected,  the $d$-foliations ${\cal F}_i$ are globally distinguishable (i.e. $N_3=d$) iff the  $d$-fold covering $W_0\to M_0$ is   trivial (i.e. $N_2=d$).
 
 On the other hand, $N_2\geq N_1$, and we can prove in fact that $N_2=N_1$ if $W$ is compact and quasi-smooth. (see [CL1] in the case $n=2$).

  \section{Global first order homogeneous polynomial partial differential equations}  
  
  In view of the the next section where we shall study the  CI-webs whose solutions satisfy to $n-1$ distinct global first order homogeneous polynomial partial differential equations, we want first to precise in this section what is a global   first order homogeneous polynomial PDE (sometimes, we shall say in short PDE in the sequel of this paper). 
  
  Above  some open set $U$ of {\bb C}$^n$, a local   first order polynomial PDE homogeneous  of degree $d$ is a differential operator  $D$ which is a holomorphic section 
   of $S^d(TU)$ (the $d^{th}$-symetric power of the tangent space $TU$),  which can always be written 
  $$ D=\sum_{|I|= d} A_I(x)\ \partial_I  ,$$        where    $I=(i_1,\cdots,i_{n},)$  denotes a multi-index of non-negative integers,  $|I|=i_1+\cdots+i_{n }$, the coefficients  $A_{I}  $   are    holomorphic functions of $x $, and   $\partial_I $ denotes the  symetric product $\prod_\lambda \Bigl({{\partial  }\over{ \partial x_  \lambda }}\Bigr)^{i_ \lambda}$.
 
 \n Solutions of this PDE are hypersurfaces $ \Sigma$ of equation $f(x_1,\cdots,x_n)=0$ in $U$, such that  $Df=0$, i.e.:   $$  \sum_{|I|= d} A_I(x). \partial_I f=0.$$ 
  For avoiding extra-solutions   defined by the vanishing   of a common factor to all of the coefficients $A_I$, 
 we shall also assume that, at any point of $U$,  the germs  of the $A_I$'s are relatively prime.  
  
 In particular, if $\Sigma $ is a graph-hypersurface of equation $\varphi(x_1,\cdots,x_{n-1})-x_n=0$, $\varphi$ is solution of the equation $$F \Bigl( x_1,\cdots,x_{n-1},\varphi(x_1,\cdots,x_{n-1}), {{\partial \varphi }\over{{\partial x_1}}} ,\ldots, {{\partial  \varphi}\over{{\partial x_{n-1}}}} \Bigr)=0,      $$
 where  $F  :\widetilde U\to \hbox{\bb C}$ is the    holomorphic function   
    of $2n-1$ variables    $ ( x_1,\ldots,x_{n-1},x_n, p_1 ,\ldots, p_{n-1})$ defined as 
    $$F (x,p)=\sum_{|I|= d}(-i)^{i_n}\ A_I(x)\ p^{\widehat I} ,$$    { where } $\widehat I=(i_1,\cdots,i_{n-1})$ will be written shortly $\widehat I= I\setminus \{i_n\}$, and $p^ {\widehat I}=(p_1)^{i_1}\ldots (p_{n-1})^{i_{n-1}} .$ 
    
    Setting $$\begin{matrix} \widetilde{F}(x, u_1,\cdots,u_{n-1},u_n)&=&(-u_n)^dF(x,-{{u_1}\over {u_n}},\cdots,-{{u_{n-1}}\over {u_n}}),\\ &&\\ &=& \sum_{|I|= d} A_I(x). u^I , 
  \hbox{ where }u^I =\prod _\lambda(u_\lambda)^{i_\lambda} ,   \\  \end{matrix}$$
    we still can write 
     $$F (x,p)=\widetilde{F}(x, p_1,\cdots,p_{n-1},-1).$$

    Let $S$ be the hypersurface of equation $F (x,p)=0$ in $\widetilde U$. 
 We  shall assume that this equation  is reduced. 
Hence, if $G(x,p)=0$ is another equation of $ {S}$, obtained from another differential operator $D_1$, there   is a unique  unit  $\rho :U\to ${\bb C}$^*$, such that  $D_1=\rho\  D $. The two equations $Df=0$ and $D_1f=0$ having the same solutions, we are mainly interested by the surface $S$.  

  If we get  $G$ and  $\widetilde{G}$ from $G$ associated to $D_1=\rho\ D$ when using  other  local coordinates $(x')$ on another open set $U'$, we must have 
  $\widetilde{G}(x',u')= \rho(x)\ \widetilde{F}(x,u)$ above $U\cap U'$.   Since $\frac{u'_n}{u_n}=\frac{  \Delta_{n-1}^ {\cal T}(x,x')}{  \Delta_{n}^ {M}(x,x')}$ after proposition 2-5, we get:
  
  \begin{lemma}  There exists some unit $\rho:U\cap U'\to ${\bb C}$^*$, uniquely defined,  such that 
 $$ G ( x' , p')=\biggl(\frac{  \Delta_{n-1}^ {\cal T}(x,x')}{  \Delta_{n}^ {M}(x,x')}   \biggr)^{- d}. \ \rho(x)\ .\ F(x,p).$$
    \end{lemma}
     
    We shall  therefore define  a global first order homogeneous polynomial     PDE of degree $d$  on a holomorphic $n$-dimensional manifold $M$ by the data of a family  $(U_a ,D_a)$ such that:

 - the  $U_a $'s are  open sets in $M$ making a covering of $M$, 
 
 -  for any $a$,    $D_a  $  is a local homogeneous polynomial PDE's of degree $d$ over $U_a$,         

 - for any pair $( a,b)$ such that $U_a\cap U_b $ be not empty, there exists a (necessarily unique) unit   $ \rho_{ ab}:U_ a \cap U_b \to  \hbox{\bb C}^*$ such that   $D_b = \rho_{ab }\ D_a$.

\n Because of their uniqueness, the functions $ \rho_{ab }$  satisfy to the cocycle condition, and are therefore    a system of transition  functions for a holomorphic line-bundle  $E $ over $M$. Hence,  
the sections  $D_a$ glue together to define a holomorphic section $D $  of the vector bundle $E\otimes S^d(TM)$  over $ M$.

Let $( ^ax$ be a system of holomorphic local coordinates on $U_a$. From the previous lemma, we deduce that the corresponding functions $F_a$ ad $F_b$ are related on $\widetilde U_a \cap \widetilde U_b$ by the formula 
$$ F_b =\biggl(\frac{  \Delta_{n-1}^ {\cal T}(^ax,^bx)}{  \Delta_{n}^ {M}(^ax,^bx)}   \biggr)^{- d}. \ \rho_{ab } \ .\ F_a.$$
Hence, the functions $F_a $ glue together, defining  
   a section $s$ of the bundle $\pi^{-1}(E)\otimes  {\cal L} ^ d$ over $\widetilde M$.  
   This justifies the following
 \begin{definition} A holomorphic line-bundle $E$ on $M$ and an integer $d\geq 1$ being given, a global   first order    polynomial PDE homogeneous of degree $d$ and     type $E$   on the  holomorphic $n$-dimensional manifold $M$  is a  
 holomorphic section $D$ of $E\otimes S^d(TM)$,
   the corresponding local equations  $F=0$ being   reduced, and  having coefficients $A_I(x)$ not depending\footnote{Notice   that    any section of $\pi^{-1}(E)\otimes {\cal L}^ d$ on $\widetilde  U$ may be written locally $ \sum_{ |I|= d} \ A_I(x,p)\ p^I\ (\sigma_{\cal L})^d\otimes \pi^{-1}(\sigma_E)$, once given a local trivialization $\sigma_E$ of $E$ and local coordinates, with the convention $p_n=-1$. Their coefficients $A_I $ depend in general on $p$, but here they don't.} on the coordinates $p$,  and having  germs at any point   relatively prime $($these conditions not depending on the local holomorphic coordinates on $M$,  and on a local holomorphic trivialization $\sigma_E$ of $E )$.  The integer $d$ is also called the \emph{weight} of the PDE.
 \end{definition}

 \begin{remark}  
  It is equivalent to give $S$, or $($up to multiplication by a global unit on $M$, i.e. a constant if $M$ is compact$)$ the section $s$, or the differential operator $D$ with coefficients in $E$. 
\end{remark}
  
\subsection{Linearizability and co-critical set of a polynomial     PDE}

Let $(D, s ,S)$ be a global polynomial PDE on $M$ as above.  

\begin{definition}  This PDE is said to be \emph{ linearizable} if it satisfies to the following assumption :

For any point $\widetilde m_0\in S$ above  a point $m_0\in M$, there exists local coordinates $x=(x_{  1},\cdots, x_{  n})$ near $m_0$, such that if  $\widetilde m_0$ has coordinates $(x_0,p_0)=(x^0_{  1},\cdots, x^0_{  n},p^0_{  1},\cdots, p^0_{  n-1})$, then the hypersurface of   equation $$x_n-x^0_{  n}=\sum_\alpha p^0_{ \alpha}(x_ \alpha-x^0_{\alpha })$$ is a solution of the PDE.  

 \end{definition} 
 \begin{definition}   We call \emph{co-critical set} of the previous PDE the analytical set $\Delta_S$ of points $\widetilde{m}\in S$
 where either $S$ is singular, or ${\cal T}_{\widetilde{m}}$ $($seen as a subspace of  ${  T}_{\widetilde{m}}\widetilde{M})$ is included into ${  T}_{\widetilde{m}}S$.
 Locally, $\Delta_S$ is defined by the equations $F=0$ and $  \frac{\partial F}{\partial x_\alpha}+p_\alpha\ \frac{\partial F}{\partial x_n}= 0  \hbox{ for any  }\alpha.$
 \end{definition} 
 \begin{proposition}  Assume a global polynomial PDE  to be linearizable. Then
its co-critical set $\Delta_S$ is all of $S$.
This  condition   may be written locally, relatively to any    system of local coordinates, in the following way$:$ 

    \n For any  $\alpha$, there exists a $($local$)$ holomorphic function $C_\alpha$ such that  $$ \frac{\partial F}{\partial x_\alpha}+p_\alpha\ \frac{\partial F}{\partial x_n}=C_\alpha\  F  .$$
 \end{proposition}
 \n Proof: Assume the web on $M$ to be locally defined by the equation 
 $F (x,p)=0$  and let \break $(x_0,p_0)=(x^0_{  1},\cdots, x^0_{  n},p^0_{  1},\cdots, p^0_{  n-1})$ be  such that  $F (x_o,p_o)=0$. If  the  hypersurface of equation 
 $$x_n-x^0_{  n}=\sum_\alpha p^0_{ \alpha}(x_ \alpha-x^0_{\alpha })$$ is a solution of the PDE, then :
 $$F\Bigl(x_1,\cdots,x_{n-1}\ ,\ x^0_{  n}+\sum_\alpha p^0_{ \alpha}(x_ \alpha-x^0_{\alpha })\ ,\ p^0_{  1},\cdots, p^0_{  n-1} \Bigr)\equiv 0.$$
By derivation of this identity,  with respect to  each $x_\alpha $, we get : 
  $\frac{\partial F}{\partial x_\alpha}+p_\alpha\ \frac{\partial F}{\partial x_n}=0 \hbox{  at the point }(x_0,p_0).$ 
 Hence $$\frac{\partial F}{\partial x_\alpha}+p_\alpha\ \frac{\partial F}{\partial x_n}\equiv 0 \hbox{ on $S$ for any $\alpha$.}$$
 But this means precisely that the vector fields   $X_\alpha=\frac{\partial  }{\partial x_\alpha}+p_\alpha\ \frac{\partial  }{\partial x_n}$, which generate ${\cal T}$ locally  are tangent to $S$ at any point of $S$. Since this property has an intrinsic meaning, not depending on the local coordinates, the same equations will be satisfied when written relatively to any other   system of local coordinates.

 \rightline{QED}
 
 \begin{proposition}Let $(D, s ,S)$ be a global polynomial PDE on $M$ as above. 
 Let $\na$ be any connection on  $\pi^{-1}(E)\otimes  {\cal L} ^ d$, and $Ds:{\cal T}\to \pi^{-1}(E)\otimes  {\cal L} ^ d$ the restriction of the covariant derivative \break $\na s:T\widetilde M\to \pi^{-1}(E)\otimes  {\cal L} ^ d$ to the sub-bundle ${\cal T}$ of $T\widetilde M$. Then, the morphism \break $Ds\vert_S:{\cal T}\vert_S\to \pi^{-1}(E)\otimes  {\cal L} ^ d\vert_S$ does not depend on the connection $\na$, and the cocritical set of the PDE is exactly the set of points $\widetilde m\in S$ where $Ds\vert_S$ vanishes. 
 \end{proposition}
 \n Proof : The section $s$ may be written locally $F.\sigma$ for a convenient trivialization $\sigma$ of the line-bundle $\pi^{-1}(E)\otimes  {\cal L} ^ d$, hence $\na s=dF\ \sigma+F\ \na \sigma$. This formula proves already that the restriction $\na s\vert_S$ of $\na s$ to $S$ (where $F$ vanishes) does not depend on $\na$. Moreover, since ${\cal T}$ is generated by the vector fields $X_\alpha= \frac{\partial  }{\partial x_\alpha}+p_\alpha\ \frac{\partial  }{\partial x_n}$, the condition $\frac{\partial F}{\partial x_\alpha}+p_\alpha\ \frac{\partial F}{\partial x_n}=0$ for any $\alpha$ or $Ds\vert_S=0$ are equivalent. 
 
 \rightline{QED}

 
\section{Complete  intersection webs    }

These webs are those whose leaves   are  the solutions of a system  of $n-1$ global polynomial partial differential equations:
$W=\bigcap_{\alpha=1}^{n-1}S_\alpha$.

 \begin{definition}
 
  A    $d$-web $W$ on $M$ is said to be a complete intersection web $($CI-web$)$ if there exists
a family $(S_\alpha  )_ \alpha$ of  $n-1$ global  polynomial  PDE's,  of respective type $E_ \alpha$ and weight   $d_ \alpha $,    such that
$d=\prod_ \alpha d_ \alpha$,   $W=\bigcap_ {\alpha=1}^{n-1} S_ \alpha$,     the sheaf of ideals of   germs of holomorphic functions vanishing on $W$ being  generated by the germs of functions vanishing on one of the $S_ \alpha$'s.
   
 \end{definition}

 If $s_ \alpha$ is a holomorphic section of $\pi^{-1}(E_ \alpha)\otimes {\cal L}^{d_ \alpha}$, such that $S_ \alpha=(s_ \alpha)^{-1}(0)$, $s_ \alpha$ may be locally written $s_ \alpha=F_\alpha \ .\ \sigma_\alpha$, where $\sigma_\alpha$ denotes a local trivialization of $E_\alpha$ (holomorphic section nowhere zero), $F_\alpha$ denoting a holomorphic function   of the shape $F_\alpha(x,p)=
 \sum_{|I|= d_\alpha} A_I^\alpha(x) p^I$ (with the convention $p_n=-1)$), the germs of the coefficients $A_I^\alpha(x)$ being relatively primes, and the equation $F_\alpha=0$ being reduced. Thus the ideal of germs of functions vanishing on $W$ is generated  by the $F_\alpha$'s.

  \begin{definition}
 A  $d$-web $W$ on $M$ is said to be a locally complete intersection web $($LCI-web$)$ if there exists a covering of $M$ by open sets $U_\lambda$ such that every induced web $W|_{ {U_\lambda}}$ is a CI-web.
  \end{definition}
  
  \begin{remark} {\rm The property for a web to be    CI   (resp.    LCI) is stronger than the property fo the underlying analytical space $W$ to be   $CI$ (resp. $LCI$) in $\widetilde M$. For instance, given $d$ distinct points $(p_i,q_i)$ on a proper conic  in the $(p,q)$-plane,   let   ${\cal F}_i$ be the the regular foliation  on  $\hbox{ \bb C} ^3$ (with coordinates $(x,y,z)$) defined by the 1-form $dz-p_idx-q_idy$.       If $d$ is odd, the regular $d$-web   $\hbox{ \bb C} ^3$ defined by these  $d$ foliations cannot be a CI-web, not even a LCI-web, while the corresponding manifold $W$ is LCI since it has no singularity. }
 \end{remark}

  \subsection{Critical   scheme  of a CI web :}
    
 In view of the study below of the dicriticity, we shall be interested in defining $\Gamma_{_W}$ not only as an analytical set, but as a scheme. This will be done now in the case of a CI web. 
 
  If  a LCI-web $W$ is  locally defined by the $n-1$ equations $F_ \alpha(x,p)=0$, $(1\leq  \alpha\leq n-1  )$, the points of $W'$ are those at which the jacobian matrix $  {\ D\ (\ F_1\ ,\ \ldots \ldots\ ,\ F_{n-1}\ ) \over{D(x_1,\ldots,x_n,p_1, \ldots,p_{n-1})}}$ has rank $n-1$, and $W_0$ is the subset of $W'$ where the
determinant $\Delta _p$ of  $ {D(F_1\ ,\ \ldots\ ,\ F_{n-1} )\over{D
(p_1\ ,\ \ldots\ ,\ p_{n-1} )}} $ does not vanish.

If we change local coordinates and  trivializations of the $E_\alpha$'s, we get for any $ \alpha$, after lemma 4-1 : 
 $$ G_  \alpha\biggl( x' ,     - \ {{  \sum_{ \beta }  K^  1 _\beta \  p_\beta-K_n^ 1  }\over{{ \sum_{ \beta }  K_ \beta ^n    \ p_ \beta-K_n^n   }}}   ,\cdots,   - \ {{  \sum_{ \beta }  K^{n-1}_\beta \  p_\beta-K_n^ {n-1}  }\over{{ \sum_{ \beta }  K _ \beta ^n    \ p_ \beta-K_n^n   }}}  \biggr)=\biggl(\frac{  \Delta_{n-1}^ {\cal T}(x,x')}{  \Delta_{n}^ {M}(x,x')}   \biggr)^{- d_  \alpha}. \ \rho_ \alpha(x)\ .\ F_ \alpha(x,p).$$ 
 Denoting  by   $ \Delta'_{p'}$ the determinant of   $ {D(G_1\ ,\ \ldots\ ,\ G_{n-1} )\over{D
(p'_1\ ,\ \ldots\ ,\ p'_{n-1} )}} $, and by $v $ the determinant of  $ {D(p'_1\ ,\ \ldots\ ,\ p'_{n-1} )\over{D
(p_1\ ,\ \ldots\ ,\ p_{n-1} )}} $, we deduce the equality  
$$ \Delta_{p } =\Bigl(\prod_\alpha  \rho_\alpha\Bigr)^{-1}. \biggl(\frac{  \Delta_{n-1}^ {\cal T}(x,x')}{  \Delta_{n}^ {M}(x,x')}   \biggr)^{-\sum_\alpha d_\alpha} \ .\ v\ .\   \Delta'_{p'}  
+ \hbox{ terms vanishing on } W. $$ 
Hence, for a CI web, the $ \Delta_p$'s glue together above  $W$,   defining a section $s_\Gamma$ over $W$     of the bundle $$\pi^{-1}\bigl(\otimes _\alpha E_\alpha\bigr)\otimes {\cal L}^{^{\sum_\alpha d_\alpha}}\otimes \bigwedge^{n-1}{\cal V}^*, $$
the zero set of which is $\Gamma_{_W}$. 
 
While the system   $(F_\alpha=0)$ of local equations for $W$ is reduced, we cannot at all assert that the system $(F_\alpha=0 \ ,\   \Delta_p=0)$ of local equations for $\Gamma_W$ has the same property. 
\begin{definition} We call \emph{critical scheme} the scheme   defined by the vanishing of the above section $s_\Gamma$.
\end{definition}

\subsection{Dicriticity  of a CI web :}

Still assuming the CI-web $W$   locally defined by the $n-1$ equations $F_ \alpha(x,p)=0$, we first want to precise the 
natural morphism ${\cal T}  |_{W_0}\to TW_0$ defining $\widetilde {\cal F}$ on $W_0$ :
hence, we are looking for the matrix $(\!(R_\alpha^\beta)\!) $ such that, for any $\alpha$, the vector field   $\frac{\partial}{\partial x_\alpha}+p_\alpha\ \frac{\partial}{\partial x_n}+\sum_\beta R_\alpha^\beta
\frac{\partial}{\partial p_ \beta}$ be tangent to $W_0$. The $R_\alpha^\beta$'s are then given by the solution of the linear system 
$$  \sum_\beta R_\alpha^\beta
\frac{\partial F_\gamma}{\partial p_ \beta} =-\Bigl(\frac{\partial F_\gamma}{\partial x_\alpha}+p_\alpha\frac{\partial F_\gamma}{\partial x_n}\Bigr)\ ,\ \hbox{ for any }\alpha,\gamma.$$
This can still be written :
$$   {D(F_1\ ,\ \ldots\ ,\ F_{n-1} )\over{D
(p_1\ ,\ \ldots\ ,\ p_{n-1} )}}\scirc ^tR =- \Theta,$$
where $\Theta$ denotes the matrix with coefficients $\Theta_ \alpha^ \gamma=\Bigl(\frac{\partial F_\gamma}{\partial x_\alpha}+p_\alpha\frac{\partial F_\gamma}{\partial x_n}\Bigr)$, hence the solution  
 $$^tR=-L\scirc \Theta    ,  \hbox { where } L= \Delta_p.\biggl({D(F_1\ ,\ \ldots\ ,\ F_{n-1} )\over{D
(p_1\ ,\ \ldots\ ,\ p_{n-1} )}}\biggr)^{-1}, $$ the matrix  $ {D(F_1\ ,\ \ldots\ ,\ F_{n-1} )\over{D
(p_1\ ,\ \ldots\ ,\ p_{n-1} )}} $   being  invertible     on  $W_0$. 

   Assuming for the moment  $W$ to be irreducible,  the only case where it is possible to extend \break ${\cal T}  |_{W_0}\to TW_0$ as an injective morphism ${\cal T}  |_{W'}\to TW'$ is the case where
all coefficients of the matrix $L\scirc \Theta $ contain $ \Delta_p$ as a common factor, i.e. vanish on the critical scheme.
In the general case, the same can be done for any irreducible component. We thus get the 
\begin{proposition} For a LCI-web $W$ to be dicritical $($resp. dicritical on an irreducible component $C$ of $W )$, it is necessary and sufficient that the matrix $L\scirc \Theta $ vanishes on the critical scheme $($resp. on  the critical scheme of $C )$.
\end{proposition}

\begin{remark} We shall see below that, in case of a linear  web on the complex projective space {\bb P}$_n$, the matrix $\Theta $ itself vanishes on the critical scheme : such a property will be called \emph{hyper-dicriticity}.  In fact, in this case, $s_\Gamma$ can be defined on all of  $\widetilde{\hbox{\bb P} _n}$,  $\Theta $ vanishes on all of $W$,  on the condition to use only affine coordinates as local coordinates. 
\end{remark}

 \begin{definition}  
   A \emph{linearizable web} on $M$ is a web   locally isomorphic to a linear web : this means that there exists local holomorphic coordinates near any point of $M_0$, with respect to which the leaves of the web have an affine equation. 
  If there   exists local holomorphic coordinates near any point of $M $ $($not only of $M_0 )$, with respect to which the leaves of the    web  have an affine equation, we shall say that the web is \emph{totally linearizable}.

     \end{definition}  

 \begin{proposition}  
 
  \hb    Any totally linearizable  web on a manifold $M$   is dicritical. 
  
  \end{proposition} 
 Proof:
  Assume a web on $M$ to be locally defined by the equations 
 $F_\alpha(x,p)=0$, ($1\leq \alpha\leq n-1$), and let $(x_0,p_0)=(x^0_{  1},\cdots, x^0_{  n},p^0_{  1},\cdots, p^0_{  n-1})$ be  such that, for any $\alpha$,   $F_\alpha(x_o,p_o)=0$. If  the   leaves of the web have an affine equation  with respect to the local coordinates $x$, the hypersurface of equation $x_n-x^0_{  n}=\sum_\alpha p^0_{ \alpha}(x_ \alpha-x^0_{\alpha })$ is a leaf, so that, for any $\alpha$, $$F_\alpha(x_{ 1},\cdots, x_{ n-1},x^0_{ n}+ \sum_\alpha p^0_{ \alpha  }(x_ \alpha-x^0_{ \alpha }),p^0_{ 1},\cdots, p^0_{ n-1})\equiv 0.$$
By derivation of this identity with respect to each $x_\beta$, we get $\Theta_\alpha^\beta=0$ at the point $(x_0,p_0)$,
hence $\Theta\equiv 0$; a fortiori $L\scirc \Theta\equiv 0$. Since, on the critical scheme,  this equation does not depend on the local coordinates, the web is dicritical as far as the local coordinates are available on $\Gamma_W$, i.e. if the web is hyperlinearizable. 

\rightline{QED}
  
\section {Webs   on the complex projective space {\bb P}$_n$.  }

\subsection {The manifold $\widetilde {\hbox{\bb P}}_n$ of contact elements  :}

  Denote by  $(X_0,\cdots,X_n)$ the homogeneous coordinates on {\bb P}$_n$, and $(u_0,\cdots,u_n)$ the homogeneous coordinates on the dual projective space  $ {\hbox {\bb P}}'_n $   of the  projective hyperplanes in  {\bb P}$_n$ : the hyperplane of    coordinates  $(u_0,\cdots,u_n)$ is the hyperplane of equation $u_0X_0+u_1X_1+\cdots+u_nX_n=0$ in  {\bb P}$_n$.
\begin{lemma} 
   The manifold $\widetilde {\hbox{\bb P}}_n$ is naturally bi-holomorphic   to the submanifold of points\break 
    $\bigl([X_0,\cdots,X_n],[u_0,\cdots,u_n]\bigr)$ in 
    $ {\hbox {\bb P}}_n\times {\hbox {\bb P}}'_n$ such that $u_0X_0+u_1X_1+\cdots+u_nX_n=0$ : a contact element  is a pair $(m,h)$ made of a point $m\in  {\hbox {\bb P}}_n$ and of a hyperplane $h\in 
{\hbox {\bb P}}'_n$     going through this  point $( m\in h)$. According to this identification, the projection  $\pi$ becomes the restriction to  $\widetilde{{\hbox {\bb P}}_n}$ of the first projection of  $ {\hbox {\bb P}}_n\times {\hbox {\bb P}}'_n$. 
    \end{lemma}
\n Proof  : We identify obviously the pair $(m,h)$  ($m\in h$)   to the sub vector-space $  T_mh$ of $T_m${\bb P}$_n$ (which is a point in {\bb  G}$_{n-1}${\bb P}$_n$).  

\begin{remark}

\hb $(i)$ If $m$ has  affine coordinates $(x_\lambda ={X_ \lambda\over X_0}$ on the open set  $X_0\neq 0$, and $h$ has for equation \break $x_n-\sum_{\alpha =1}^{n-1}p_\alpha x_\alpha=0$ with respect to these coordinates, we identify $\bigl((x_\lambda)_\lambda,(p_\alpha)_\alpha\bigr)$ to the point \break $\bigl([X_0,\cdots,X_n],[u_0,\cdots,u_n]\bigr)$ with 
$X_0=1$, $X_\lambda=x_\lambda$, $u_0= x_n-\sum_{\alpha =1}^{n-1}p_\alpha x_\alpha$,  $u_\alpha=p_\alpha$ and $u_n=-1$. 

\hb $(ii)$ Conversely, we identify the point  $\bigl([X_0,\cdots,X_n],[u_0,\cdots,u_n]\bigr)$ in the open set $X_0. u_n\neq 0$ of $ {\hbox {\bb P}}_n\times {\hbox {\bb P}}'_n$,  such that $u_0X_0+u_1X_1+\cdots+u_nX_n=0$,    with the point of coordinates $\bigl(x_\lambda ={X_ \lambda\over X_0}, p_\alpha=-{u_ \alpha\over u_n}\bigr)$. 
\end {remark}

 The   spaces 
{\bb P}$_n$ et $ {\hbox {\bb P}}'_n$  play    the same role, so that $\widetilde {\hbox{\bb P}}_n= \widetilde {\hbox{\bb P}}'_n$   : the second  projection\break  $ \pi':\widetilde{{\hbox {\bb P}}_n}\to  {\hbox {\bb P}}'_n$ is also a fiber space   with fiber {\bb P}$_{n-1}$,  with which we can do the same constructions as with  $\pi$. Let ${\cal V}'$, ${\cal T}'$ and ${\cal L}'$ be the bundles built from $ \pi'$ in the same way as  ${\cal V}$, ${\cal T}$ and ${\cal L}$ are built from  $ \pi$. Writting $p_n=x_n -\sum_\alpha p_\alpha x_\alpha$, notice that,   the 1-form $dx_n-\sum_\alpha p_\alpha dx_\alpha$ may be written $dp_n-\sum_\alpha x_\alpha dp_\alpha$ on the open set $X_0. w_n\neq 0$ of $\widetilde {\hbox{\bb P}}_n$, with respect to the coordinates $(-x_\alpha, p_n, p_\alpha)$. The tautological contact form 
$\omega$ is then the same when $\widetilde{{\hbox {\bb P}}_n}$ is seen as the manifold of contact elements of $\widetilde{{\hbox {\bb P}}_n}$ or of $\widetilde{{\hbox {\bb P}}'_n}$. Therefore, when a $n$-dimensional subvariety $W$ of  $\widetilde{{\hbox {\bb P}}_n}$ is simutaneously a web on ${\hbox{\bb P}}_n$ and on ${\hbox{\bb P}}'_n$ (``bi-web"), the foliation $\widetilde {\cal F}$ on the regular part of $W$ is the same for both webs.

Denote  respectively by  ${\cal O} (  - 1)$ and  ${\cal O} '( -  1)$  the  tautological line bundles of  ${\hbox {\bb P}} _n$ and  ${\hbox {\bb P}}'_n$, ${\cal O} (  1)$ and  ${\cal O} '(  1)$ their  dual bundle. Set    $\ell={\pi}^{-1}\bigl( { {\cal O}}( -  1)\bigr)$ and  $\ell'={\pi'}^{-1}\bigl( {  {\cal O'}}( - 1)\bigr)$. More  generally, $\ell^{ k}$ (resp. $\ell^{- k}$ or $\breve l^k$) will denote, for any integer $k\geq 0$,the  $k^{th}$ tensor power of  $\ell$ (resp. of its    dual bundle $\breve \ell$), and the same goes  for  $(\ell')^{ k}$ and  $(\ell')^{- k}=(\breve  l')^k$.
\begin{lemma}
\hb $(i)$ The line bundle  $ \cal  L $ is  isomorphic to $\breve  \ell  \otimes \breve\ell'$, which is also the normal bundle $N_{\widetilde{\hbox {\bb P} _n}}$ of $\widetilde{\hbox {\bb P} _n}$ in $ {\hbox {\bb P} _n}\times  {\hbox {\bb P}' _n}$. 
\hb $(ii)$ The line bundle $ \pi^{-1}\Bigl( \bigwedge ^n T{\hbox {\bb P} _n}\Bigr) $   is  isomorphic to $(n+1)\breve  \ell  $. 
\hb $(iii)$ The sub-bundles  ${\cal V}$ and ${\cal T}'$ of $T\widetilde{ \hbox {\bb P}}_n$ the same goes ,  as well as ${\cal V}'$ and ${\cal T}$.     

   \end{lemma}
\noindent Proof :

On the open set $(X_0\neq 0)\cap (X_n\neq 0)$ of $ { \hbox {\bb P}}_n$, we get    $p'_\alpha={{-p_\alpha} \over {x_n-\sum _\alpha p_\alpha x_\alpha}}  $ 
by  the change of affine coordinates  $(x'_\alpha={ {x_\alpha}\over {x_n}},x'_n={
1\over {x_n}})$,  hence the formula 
$$dx _n -\sum _\alpha p _\alpha\ x _\alpha= -x_n\bigl(x_n-\sum_\alpha p_\alpha\ x_\alpha\bigr)\Bigl(dx'_n  -\sum _\alpha p'_\alpha\ x'_\alpha\Bigr).$$
Since $-x_n\bigl(x_n-\sum_\alpha p_\alpha\ x_\alpha\bigr)={{-X_n}\over{X_0}}\bigl({{\sum_\lambda u_\lambda X_\lambda}\over{u_nX_0}}\bigr)$, and since $\sum_\lambda u_\lambda X_\lambda=-u_0X_0$ on  $\widetilde {\hbox{\bb P}}_n$, we get finally $-x_n\bigl(x_n-\sum_\alpha p_\alpha\ x_\alpha\bigr)={{X_n}\over{X_0}}{{u_0}\over{u_n}}$. Thus, ${{X_n}\over{X_0}}{{u_0}\over{u_n}}$ is the transition function for ${\cal L}^*$  on  the open set\break  $(X_0.u_n\neq 0)\cap (X_n.u_0\neq 0)$; but this is also   the transition function for the bundle $\ell\otimes   \ell'$, hence part $(i)$   the lemma.

We know in general  that  the bundles $TM$ and $ \bigwedge ^n TM$ have the same Chern class $c_1$. Here   the formula  $T { \hbox {\bb P}}_n\oplus 1=(n+1){\cal O}(1)$ implies that  $c_1(T { \hbox {\bb P}}_n)$ is equal to $ c_1((n+1){\cal O}(1))$. But this Chern class is still equal to that of ${\cal O}(n+1)$. Hence, the bundles $ \bigwedge ^n T { \hbox {\bb P}}_n$ and 
${\cal O}(n+1)$ are isomorphic, since the isomorphy class of a complex line-bundle is characterized by its Chern class. This proves part $(ii)$ of the lemma. 

Using the coordinates $(x_\lambda ={X_ \lambda\over X_0}\ ,  \  p_ \alpha =-{u_ \alpha\over u_n})$ on the open subset $X_0.u_n\neq 0$ of $\widetilde{ \hbox {\bb P}}_n$,  we make the  change of coordinates : 
$$x'_\alpha=p_\alpha\ ,\ x'_n=x_n-\sum_\alpha p_\alpha\ x_\alpha\ ,\  p'_\alpha=x_\alpha.$$ Notice that $x'_n$ is also equal to $-{u_0\over u_n}$. We get the formulae:
$$\frac{\partial}{\partial p'_\alpha}=\frac{\partial}{\partial x_\alpha}+p_\alpha \frac{\partial}{\partial x_n}\hbox{ \ and \ }\frac{\partial}{\partial p _\alpha}=\frac{\partial}{\partial x'_\alpha}+p'_\alpha \frac{\partial}{\partial x'_n}.$$
Then, any vector of ${\cal V}'$ is a vector of ${\cal T}$  and conversely, when $X_0.u_n\neq 0$. Any point of $\widetilde{ \hbox {\bb P}}_n$ belonging to some open subset $X_i.u_j\neq 0$ (with $i\neq j$, $0\leq i,j\leq n$, the previous identification may be done above all of $\widetilde{ \hbox {\bb P}}_n$. We have of course the similar identification between 
 ${\cal V} $ and ${\cal T}'$, hence the part $(iii)$ of the lemma.

\subsection{Degree of a web on  {\bb P}$_n$, and multi-degree of a CI-web  :}

\subsubsection{Bi-degree of a global   PDE  on  {\bb P}$_n$  :}

Let  $H( X_0, \cdots,X_n;u_0, \cdots,u_n) $ be a  polynomial with respect to the variables  $( X_0, \cdots,X_n;u_0, \cdots,u_n)$    homogeneous of degree  $ \delta$ with respect to the  variables $X_0, \cdots,X_n$, and   homogeneous of degree  $d$ with respect to the  variables $u_0, \cdots,u_n$. We shall call  $( \delta,d)$ the  \emph {bi-degree} of  $H$. Let  $S$ be the  hyper-surface of equations $(H=0,\sum_{\rho=0}^n u_\rho\ X_\rho=0)$ in $\widetilde{{\hbox {\bb P}}_n}$.  
Any other polynomial $\overline H$ defining the same hyper-surface    $S$ must have the same  bi-degree. 

Such a hypersurface $S$ is the zero set of a holomorphic section of $\breve\ell ^\delta\otimes\breve(\ell') ^d$, or equivalently (after lemma 6-3-$(i)$)   of $\pi^{-1}(E)\otimes {\cal L}^d$, with $E={\cal O}(\delta-d)$. It defines a global  first order PDE on {\bb P}$_n$, homogeneous of degree $d$.

 The integer  $ \delta$ (resp. $d$) is in fact  the  degree  of the algebraic $(n-2)$-dimensional subvariety of the points $m\in $ {\bb P}$_n$ (resp. of the points $h\in $ {\bb P}'$_n$ such that $(m,h)$ is a formal solution of the given PDE at order 1, $h $  denoting a given  generic hyperplane  of {\bb P}$_n$ (resp. $m $  denoting a given  generic point  in {\bb P}$_n$).

By restriction to the set   $X_0.u_n\neq 0$, and  relatively to the corresponding affine coordinates \break $x_\lambda=X_\lambda/X_0, p_ \alpha=-u_ \alpha/u_n$,
 $S$ has  equation 
 $$H\bigl(1,  x_1, \cdots,x_n\ ;\ x_n-\sum_\alpha p_\alpha \ x_\alpha\ ,\ p_1, \cdots,p_{n-1},-1\bigr)=0, $$
the first member of which  
is a polynomial of degree $d$ in $p$, with polynomial coefficients 
   $a_i(x)$ of  degree $\leq  \delta+d$ in $x$.

Conversely, any global $PDE$ on    ${\hbox {\bb P}}_n$ may be defined by this procedure from a polynomial $H$. Identifying   {\bb C}$^n$ to the  affine open set $X_0\neq 0$ in  ${\hbox {\bb P}}_n$, with affine coordinates $x_\lambda={X_\lambda\over X_0} $.
 \begin{lemma}  
 \hb $(i)$ A global PDE on  ${\hbox {\bb P}}_n$ is completely defined by its  restriction to {\bb C}$^n$.
 \hb $(ii)$ For a  global $($first order, homogeneous polynomial$)$ PDE on   {\bb C}$^n$ to be the restriction of a PDE on ${\hbox {\bb P}}_n$, it is necessary and sufficient that the  coefficients
 $a_i$ be all polynomial in  $(x_1,\cdots, x_n)$.
 \end{lemma}

\n  The proof is completely similar to that given in [CL1] for  $n=2$ (notice that, when $n=2$, any homogeneous first order PDE is   automatically a web). 

\subsubsection{Algebraic PDE's on  {\bb P}$_n$  :}

 \begin{definition}  
A global    PDE $(S,D,s)$, homogeneous of degree $d$ on {\bb P}$_n$,  is said to be \emph{algebraic}, if there exists an algebraic hypersurface $\overline{S}$ of degree $d$ in {\bb P}'$_n$ $($necessarily unique$)$, such that $S$ be the set of all points $(m,h)$ with $m\in h$ and $h\in 
\overline{S}$. 
 \end{definition}  
 
  \begin{lemma}  
The algebraic polynomial PDE $S$'s, homogeneous of degree $d $,  are the PDE's of bi-degree $(0,d)$. 
\end{lemma}
  Proof :  
Let $\Phi(u_0,u-1,\cdots,u_n)=0$ be the  equation of an algebraic hypersurface of degee $d$ in {\bb P}'$_n$ (well defined up to multiplication by a non-zero scalar). Then, the PDE $S$ defined by the polynomial  $H( X_0, \cdots,X_n;u_0, \cdots,u_n) =\Phi(u_0,u-1,\cdots,u_n)$    is algebraic, and the map so defined is an obvious bijection  onto the set of global PDE's   of bi-degree $(0,d)$. [Observe that the polynomial $H$ of such a PDE of bi-degree $(0,d)$ is well defined, since $H( X_0, \cdots,X_n;u_0, \cdots,u_n).\sum_\rho u_\rho\ X_\rho$ cannot have bi-degree $(0,d)$ if $K\neq 0$. 

\rightline {QED}

\subsubsection{Multi-degree of a CI-web  on  {\bb P}$_n$  :}

A CI-web  $W$  is the complete intersection  of $n-1$ global PDE's $S_\alpha$ $(1\leq \alpha\leq n-1)$. Let $(\delta_\alpha,n_\alpha)$ be the bi-degree of $S_\alpha$ for such a CI-web. 
\begin{definition} The family $(\delta_\alpha)_\alpha$ is called the {\emph multi-degree} of the CI-web, and the number $\delta=\prod_\alpha\delta_\alpha$ its \emph{degree}.
\end{definition}

  \subsection {Linear  and algebraic  webs: }  
  
  \begin{definition}  \hb
  $(i)$ A linear web on  an open set of  {\bb P}$_n$ is a web   all leaves of which  are pieces of   hyperplanes.
  \n $(ii)$ A  algebraic $d$-web on {\bb P}$_n$ is a web whose leaves are the hyperplanes belonging to some algebraic curve of degree $d$ in {\bb P}'$_n$.
 \end{definition}

\begin{theorem} The algebraic webs on {\bb P}$_n$ are the linear webs globally defined on all of {\bb P}$_n$ 
\end{theorem}
Proof : If a web $W$ is linear, and if a point $(m_0, h_0)$ belongs to $W$, then all points $(m,h_0)$ such that $m$ belongs to $h_0$ are still in $W$. Therefore, the web is completely defined by the projection $\overline{W}=\pi'(W)$.
It is then sufficient to prove that $\overline{W}$ is an algebraic set, because it will be  then automatically one-dimensional since $W$ has dimension $n$. 

It is then possible to generate the ideal of   $W$ by functions on  $\widetilde{\hbox{\bb P}}_n$ going  to the quotient modulo $\pi'$ ; we get therefore  analytical functions on {\bb P}'$_n$ defining $\overline{W}$. Hence, $\overline{W}$ is analytic, and conseqently algebraic\footnote{We have already proved this theorem when  $n=2$ by another method (see [CL1]). The principle of the method used here  has been suggested to us by L. Pirio.}. 
\rightline{QED}

  \begin{remark} Because of the duality between the curves in {\bb P}$'_n$ and the developable hypersurfaces in {\bb P}$  _n$,  the hyperplanes of an algebraic web are also the hyperplanes tangent to some algebraic developable hypersurface $\cal C$ in {\bb P}$_n$.  This means that  there exists an algebraic curve $\gamma$ (i.e. analytical set of pure dimension 1),   a holomorphic  fiber-bundle  $\widehat {\cal C}\to \gamma$ with  base the analytical set $\gamma$ and fiber {\bb P}$_{n-2}$, and an immersion  $\Phi :\widehat {\cal C}\to ${\bb P}$_n$ of the total space $\widehat {\cal C}$ of this bundle into {\bb P}$_n$, 
  
  -   whose image 
  $\Phi(\widehat {\cal C})$ is ${\cal C}$, 
  
  - whose restriction to any fiber of $\widehat {\cal C}$ is a bi-holomorphism onto some  $(n-2)$ - dimensional sub projective space of {\bb P}$_n$, 
  
  - and such that the tangent hyperplane to ${\cal C}$ at some regular   point $\Phi(a)$ of $\cal C$ depends only on the fiber of $\widehat {\cal C}$ to which $a$ belongs. 
  \end{remark}
  
  \begin{proposition}  
 The algebraic CI-webs on {\bb P}$_n$ are the webs $W=\cap_\alpha S_\alpha$ of multi-degree $(\delta_\alpha)_\alpha$ with all $\delta_\alpha$ equal to 0.
 \end{proposition}
  Proof: It is clear that, if $W=\cap _\alpha S_\alpha$ such that all $ \delta_\alpha$'s are zero, then $W$ is algebraic 
  (it corresponds to the curve in {\bb P}$'_n$ defined by the intersections of the $ \overline{S}_\alpha$'s. Conversely, assume that  $W=\cap _\alpha S_\alpha$ is algebraic. Assume that there  exists at least   one of the   $S_{\alpha }$'s  (let us say $S_{\alpha_0}$ which is not algebraic, and let $(m_0, h_0)$ be a point of $W$.   Since  $\delta_{\alpha_0}\geq 1$, there exists $m\in ${\bb P}$_n$, such that $(m , h_0)$ does not belong to $S_{\alpha_0}$, and a fortiori does not belong to $W$. But there is a contradiction with the algebraicity of $W$ (since $(m_0, h_0)$ belongs to $W$, all points $(m,h_0)$ should belong to $W$).

\subsection{Cohomology of  $\widetilde{\hbox {\bb P}_n}$ : }

 Let us denote  respectively by    $\xi= c_1(\breve\ell  ) $ and  $\xi'= c_1 (\breve\ell')    $   the Chern classes in $ H^2(\widetilde{\hbox {\bb
P}_n},\hbox {\bb
 Z}) $ of  the bundles $\pi^{-1}\bigl( {\cal O}(1)  \bigr)$ and  $   (\pi')^{-1}\bigl(   {\cal O}'(1)  \bigr)$.  
 \begin{theorem}  The cohomology algebra of $\widetilde{\hbox {\bb P}_n}$  is the quotient of the free algebra $\hbox{\bb Z}[\xi,\xi']$ by the 3 relations $\xi^{n+1}=0$, $(\xi')^{n+1}=0$ and $\sum_{j=0}^n(-1)^j\xi^j(\xi')^{n-j}=0\  :$ 
   $$H^*(\widetilde{\hbox {\bb
P}_n},\hbox {\bb
 Z})=\hbox{\bb Z}[\xi,\xi']/\bigl(\xi^{n+1},\ (\xi')^{n+1},\ \sum_{i=0}^n(-1)^i\xi^i(\xi')^{n-i}  \bigr).$$
 \end{theorem} 
 \n Proof : Since the base   {\bb P}$_n$ and the fiber {\bb P}$_{n-1}$ of the fibration $\pi$ have only even cohomology, the spectral sequence of this fibration   collapses, so that 
$H^*(\widetilde{\hbox {\bb
P}_n},\hbox {\bb
 Z})$ is isomorphic -as a graded vector space- to the $E_2$ - term   $\Bigl(\hbox{\bb Z}[\xi]/ \xi^{n+1}\Bigr)\otimes \Bigl(\hbox{\bb Z}[ \bar \eta]/\bar \eta^{n}\Bigr)$, where $\bar \eta$ denotes the generator in $H^2( {\hbox {\bb
P}_{n-1}},\hbox {\bb
 Z}) $ of the cohomology  algebra  of the fiber.  
 
 Since ${\cal L}^*$  is the tautological line bundle of $  \widetilde{\hbox {\bb P}_n}$=
{\bb P}$( T^*M)$, the Chern class $\eta=c_1(\breve {\cal L}^*)$ induces on every fiber of $\pi$  the generator of the cohomology of this fiber. Therefore, 
   $H^*(\widetilde{\hbox {\bb
P}_n},\hbox {\bb
 Z})$ is generated by $\xi$ and $\eta$ as an algebra, or also by $\xi$ and $\xi'$  since we can deduce the equality $\eta=-(\xi+\xi')$
  from the isomorphism $ {\cal  L }^*\cong \ell   \otimes \ell'$ : hence $H^*(\widetilde{\hbox {\bb
P}_n},\hbox {\bb
 Z})$ is a quotient of the free algebra $\hbox{\bb Z}[\xi,\xi']$. 
 
 The relations $\xi^{n+1}=0$ and $(\xi')^{n+1}=0$ are obvious for dimensional reasons. 
  
 \begin{lemma}
  The  Chern class $c_j( {\cal T})$ is given by the formula : 
 $$c_j( {\cal T})= \sum_{i=0}^j(-1)^i\begin{pmatrix} n+1\\j-i\\  \end{pmatrix}\xi^{j -i}(\xi+\xi')^i.$$ 
 \end{lemma}
 \n Proof :  From the identity on the total Chern classes 
 $ c ( {\cal T})(1+ \xi+\xi')=(1+  \xi )^{n+1}$  induced by the exact sequence $(*)$ of vector bundles, we can compute $c_j( {\cal T})$ by induction on $j$ with the formula 
$c_j( {\cal T})= \begin{pmatrix} n+1\\j\\  \end{pmatrix}\xi^j -(\xi+\xi')c_{j-1}( {\cal T})  $.

 In particular, since $c_n( {\cal T})=0$, we get the relation 
  $\sum_{i=0}^n(-1)^i\begin{pmatrix} n+1\\n-i\\  \end{pmatrix}\xi^{n -i}(\xi+\xi')^i=0,$ 
  which can still be  formally written   ${{\Bigl( \bigl((\xi+\xi')-\xi\bigr)^{n+1}+\xi^{n+1}\Bigr)}/ {\bigl(\xi+\xi'\bigr)}}=0$, i.e. 
  $$\sum_{1=0}^n(-1)^i\xi^i(\xi')^{n-i}=0.    $$  
  
  Since $\hbox{\bb Z}[\xi,\xi']/\bigl(\xi^{n+1},\ (\xi')^{n+1},\ \sum_{i=0}^n(-1)^i\xi^i(\xi')^{n-i}  \bigr)$ and  $E_2$ are isomorphic as graded vector spaces, we get the formula of the theorem.

\rightline{QED} 
 
\begin{theorem}  Any CI-web     $W$ of codimension one on 
{\bb P}$_n$ has a non-empty critical set $\Gamma_W$. Equivalently, such a web has a non-empty caustic. 
\end{theorem}

  \n Proof : 

If  $\Gamma_W$ was empty, $W$ would be in particular non-singular, the foliation $\widetilde{\cal F}$ defined by the morphism $\omega_W:TW\to {\cal L}$ would have no singularity, so that  $(c_1)^2( {\cal L})\frown [W]$ would be zero after the Bott vanishing theorem.

 Since the CI web $W$ is the zero set of a section of the bundle $\widetilde{N}_W=\bigoplus _\alpha \Bigl( \breve \ell ^{ \delta_\alpha}\otimes (\breve \ell')^{d_\alpha}\Bigr)$ above $\widetilde{\hbox {\bb P}}_n $, the fundamental class $[W]$ is the Poincar\é dual of the 
   Chern class $c_{n-1}(\widetilde{N}_W)= \prod_\alpha(\delta_\alpha\xi+d_\alpha\xi')|_{_W}$. 
   Hence, $$(c_1)^2( {\cal L})\frown  [W]=\Bigl((c_1)^2( {\cal L})\smile
 c_{n-1}(  \widetilde  N_W     )\Bigr)\frown \bigl[\widetilde{\hbox {\bb P}}_n \bigr] \hbox { in } H_{2(n-2)}(\widetilde{\hbox {\bb P}}_n).$$
  Since $\widetilde{\hbox {\bb P}}_n $ is the zero set of a section of $\cal L$ on $\hbox{\bb P}_n\times \hbox{\bb P}'_n$, we still have 
 $$(c_1)^2( {\cal L})\frown  [W]=\Bigl((c_1)^3( {\cal L})\smile
 c_{n-1}(  \widetilde  N_W     )\Bigr)\frown \bigl[\hbox{\bb P}_n\times \hbox{\bb P}'_n  \bigr] \hbox { in } H_{2(n-2)}(\hbox{\bb P}_n\times \hbox{\bb P}'_n). $$

Denote  $\bar \xi$ and $\bar \xi'$   the Chern classes of  $ {\cal O}(1) $
 and $ {\cal O}'(1) $  in $H^2(\hbox{\bb P}_n)$ and $H^2(\hbox{\bb P}'_n)$ respectively, or in $H^*(\hbox{\bb P}_n\times \hbox{\bb P}'_n)=\Bigl(\hbox{\bb Z}[\bar\xi]/\bar \xi^{n+1}\Bigr)\otimes\Bigl( \hbox{\bb Z}[\bar\xi']/(\bar \xi')^{n+1}\Bigr)$. 

But
$ (c_1)^3( {\cal L})\smile
 c_{n-1}(  \widetilde  N_W     ) =(\bar \xi+\bar \xi')^3\prod_{\alpha =1}^{n-1}(\delta_\alpha \bar\xi+d_\alpha \bar\xi')$ may be written 
\medskip
$\sum_{i=2}^{n}a_i\  \bar\xi^i( \bar\xi')^{n+2-i}$ with coefficients $a_i$ all strictly positive. Hence  $(c_1)^2( {\cal L})\frown [W]$ may not vanish.

\rightline{QED}

\section{Dicriticity  and algebraicity of global CI-webs}

While all results of this section could be given for   webs   which are not necessarily 
CI or LCI, we shall restrict ourselves to CI-webs for simplicity.

 \begin{proposition}

    Any linear web on an open set of {\bb P}$_n$    $($in particular  any algebraic web on {\bb P}$_n)$ is totally linearizable, hence dicritical.
   
   \end{proposition} 
   Proof:
If a web on {\bb P}$_n$ is linear, the equation of the leaves are affine with respect to any system of affine coordinates. Since there are such affine coordinates near any point of {\bb P}$_n$ , and in particular near any point of the critical set, this proves that a linear web is hyperlinearizable, hence part (ii). 

\rightline{QED}

 For   webs   which are globally defined on  {\bb P}$_n$, we have therefore the implications
 $$algebraic \Longleftrightarrow linear \Longrightarrow dicritical. $$
 We shall prove now conversely, at least in the case of  webs whose each irreducible component is CI, that in fact 
    $$(\hbox{\it dicritical\ $+$\ quasi-smooth})\Longrightarrow algebraic.$$

    \begin{theorem}   For $d\geq 3$, any  $d$-web on {\bb P}$_n$   which is quasi-smooth and dicritical, and whose  each irreducible component is CI, is algebraic.
   
   \end{theorem}
\noindent Proof : Since quasi-smoothness  and dicriticity are defined on each irreducible component, we may assume   $W$ irreducible and smooth. Dicriticity means then that the foliation $\widetilde{\cal F}$ has no singularity. After Bott
  $(c_1)^2 \bigl(N(\widetilde{\cal F} )\bigr)$ must vanish. 
  \n The total Chern class $c  \bigl(N(\widetilde{\cal F} )\bigr)$ is equal to $c(W).c({\cal T}|_W) ^{-1}$.
  
 \n After the exact sequence $(*)$ of section 2, $c({\cal T})=\bigl(\pi^{ *}c (\hbox{\bb P}_n)\bigr).c({\cal L})^{-1} $.
 
  \n Assume the CI-web to be defined by the $n-1$ global PDE's $S_\alpha$ of bi-degree $(\delta_\alpha\geq 0, d_\alpha\geq 1)$, with $d=\prod_\alpha   d_\alpha $ and $ \delta=\prod_\alpha  \delta_\alpha $. Then
  $$c(W)=c\bigl(\widetilde{\hbox{\bb P}_n})|_W 
       . \Bigl( \prod _\alpha \bigl(1+{\delta_\alpha}\ \xi +{d_\alpha\ \xi' }\bigr)|_{_W}\Bigr)^{-1}.$$
         Since $\widetilde{\hbox{\bb P}_n} $ has ${\cal L}$ for normal bundle in $ {\hbox{\bb P}_n} \times  {\hbox{\bb P}'_n} $, and since $(\pi')^{-1} (T {\hbox{\bb P}'_n})\oplus 1=(n+1)\breve \ell$, we get finally : 
   $$c  \bigl(N(\widetilde{\cal F} )\bigr)=c(\pi'^{-1} T {\hbox{\bb P}'_n})|_{_W}\smile \Bigl( \prod _\alpha \bigl(1+{\delta_\alpha}\ \xi +{d_\alpha\ \xi' }\bigr)|_{_W}\Bigr)^{-1},$$
   and in particular 
   $$c _1 \bigl(N(\widetilde{\cal F} )\bigr)=\Bigl((n+1)\xi'-\sum_\alpha (\delta_\alpha \xi+d_\alpha \xi')\Bigr)|_{_W}.  $$
   Let us write $\overline{\delta} =\sum_\alpha  \delta_\alpha $ and $\overline{d} =\sum_\alpha d_\alpha $ . 
   We thus get 
   $$(c_1)^2 \bigl(N(\widetilde{\cal F} )\bigr) =\Bigl[(\overline{\delta})^2\xi^2
   + (n+1-\overline d)^2(\xi')^2- 2\overline{\delta}(n+1-\overline{d})\xi\xi'\Bigr]|_{_W}.$$
   Since the CI web $W$ is the zero set of a section of the bundle $\widetilde{N}_W=\bigoplus _\alpha \Bigl( \breve \ell ^{ \delta_\alpha}\otimes (\breve \ell')^{d_\alpha}\Bigr)$, the fundamental class $[W]$ is the Poincar\é dual of the 
   Chern class $c_{n-1}(\widetilde{N}_W)= \prod_\alpha(\delta_\alpha\xi+d_\alpha\xi')|_{_W}$, so that 
   $$(c_1)^2 \bigl(N(\widetilde{\cal F} )\bigr)  \frown [W] 
  =\biggl(\Bigl[(\overline{\delta})^2\xi^2
   + (n+1-\overline d)^2(\xi')^2- 2\overline{\delta}(n+1-\overline{d})\xi\xi'\Bigr] \smile  \prod_\alpha(\delta_\alpha\xi+d_\alpha\xi') \biggr)\frown  [\widetilde{\hbox{\bb P}}_n]  $$
    { in } $H^{2n+2}(\widetilde{\hbox{\bb P}_n})$. 
     If $(c_1)^2 \bigl(N(\widetilde{\cal F} )\bigr)=0$, 
    then  $$ \Bigl[(\overline{\delta})^2\xi^2
   + (n+1-\overline d)^2(\xi')^2- 2\overline{\delta}(n+1-\overline{d})\xi\xi'\Bigr]  \smile \Bigl( \prod_\alpha(\delta_\alpha\xi+d_\alpha\xi')\Bigr) \smile \xi^{n-2}    $$ 
    must vanish in $H^{4n-2}(\widetilde{\hbox{\bb P}}_n)$. 
     Using the relations $ \xi^{n+1}=0$, $( \xi')^{n+1}=0$, and  $ \xi^{n }( \xi')^{n -1}= \xi^{n-1}( \xi')^{n } $,    we deduce that the  number  
   $${\cal N}=\overline{\delta}^2-2\overline{\delta} (n+1-\overline{ d})(1+\sigma_1)+(n+1-\overline{ d})^2( \sigma_1+\sigma_2) $$ 
  must vanish, where $\sigma_1$ (resp. $\sigma_2$) denotes the first (resp the second) elementary symmetric function 
  $\sigma_1=\sum _\alpha \delta_\alpha/d_\alpha$ (resp. $\sigma_2=\sum _{\alpha \beta} \delta_\alpha \delta_\beta/d_\alpha d_\beta$) of the numbers  $\delta_\alpha/d_\alpha$.

  Since every  $d_\alpha$ is at least equal to 1, and since $  \prod_\alpha d_\alpha \ (=d)$ is at least 3,   $ n+1-\overline{d}  $ is always non-positive. Thus $ \cal N$  is   non-negative, and may be zero only if all $\delta_\alpha$'s vanish.  This is  what we want. 
  
  \rightline{QED}
 
 \n {\bf References}
 
  \n [B] W. Blaschke, {\it Einf\"uhrung in die Geometrie der Webe}, Birkha\"user, Basel, 1955. 
  
 \n [C] S.S.Chern, {\it Abz\"ahlungen f\"ur Gewebe}, Abh. Hamburg 11, 1936, 163-170.

 \n [CL1]   V. Cavalier et D. Lehmann, {\it Introduction \à l'\étude globale des tissus sur une surface holomorphe}, Annales de l'Institut Fourier 57-4 (2007),  1095-1133. 
 
\n  [CL2]   V. Cavalier et D. Lehmann, {\it Regular holomorphic webs of codimension one}, arXiv : math. \break DS 0703596, v1, 20 march 2007. 

\n [H]  A. H\énaut, {\it On planar web geometry through abelian relations and connections}, Annals of Maths 159 (2004), 425-445. 
 
 \n [Y]   J.N.A.  Yartey, {\it Number of singularities of a generic web on the complex projective plane}, J. of Dynamical and Control Systems 11-2, april 2005, 281-296.
 
  \bigskip
\n Vincent Cavalier, Daniel Lehmann\hb D\'epartement des Sciences Math\'ematiques, CP 051,
Universit\'e de Montpellier II \hb Place Eug\`ene Bataillon, F-34095 Montpellier
Cedex 5, France\hb   cavalier@math.univ-montp2.fr;\ \ \ lehmann@math.univ-montp2.fr
 
  \end{document}